\documentclass[letterpaper, 12pt]{article}
\usepackage{amsmath}
\usepackage{amssymb}
\usepackage{amsfonts}
\usepackage{mathrsfs}
\usepackage{bbm}
\usepackage{color}
\usepackage{graphicx}

\newcommand{\ecal}{\mathcal{E}}

\newcommand{\gcal}{\mathcal{G}}

\newcommand{\scal}{\mathcal{S}}

\newcommand{\vcal}{\mathcal{V}}

\newcommand{\real}{\mathbb{R}}
\newcommand{\intgr}{\mathbb{Z}}

\newcommand{\ind}{\mathbbm{1}}

\newcommand{\bdm}{\begin{displaymath}}
\newcommand{\edm}{\end{displaymath}}
\newcommand{\bea}{\begin{eqnarray*}}
\newcommand{\eea}{\end{eqnarray*}}
\newcommand{\bean}{\begin{eqnarray}}
\newcommand{\eean}{\end{eqnarray}}

\newcommand{\bfx}{\mathbf{x}}

\newtheorem{proposition}{Proposition}

\newtheorem{definition}{Definition}

\newtheorem{remark}{Remark}

\author{{\Large S\'ebastien Roch}\\
Department of Statistics\\
University of California, Berkeley\\
Berkeley, CA 94720--3860\\
\texttt{sroch@stat.berkeley.edu}}
\title{Bounding Fastest Mixing}

\begin{document}

\maketitle

\begin{abstract}
In a series of recent works, Boyd, Diaconis, and their co-authors have introduced a semidefinite programming approach for computing the fastest mixing Markov chain on a graph of allowed transitions, given a target stationary distribution. In this paper, we show that standard mixing-time analysis techniques---variational characterizations, conductance, canonical paths---can be used to give simple, nontrivial lower and upper bounds on the fastest mixing time. To test the applicability of this idea, we consider several detailed examples including the Glauber dynamics of the Ising model---and get sharp bounds. 
\end{abstract}

\bigskip


\section{Introduction}

Sampling from a complex collection of objects is a basic procedure 
in physics, statistics and computer science. A widely used
technique, known as Markov chain Monte Carlo (MCMC), 
consists in designing a Markov chain on the set to be
sampled such that the law of the chain
converges to the desired distribution. The chain 
is run long enough for a sample to be picked
from a good approximation of the stationary distribution. 
The time one has to
wait in order for this approximation to be satisfactory 
is known as the mixing time. In practice, it is crucial
that this parameter be small. See
e.g.~\cite{J03} for a survey of theoretical results on MCMC.

One way to picture a Markov chain (MC) on a combinatorial structure
is to think of the states as nodes and of the transitions
as edges. For a chain to be implementable, the neighbourhood 
structure surrounding each node must be relatively simple. Under
this constraint, one has to choose a set of allowed transitions
that is most likely to produce fast convergence. This
is usually done in a heuristic manner.

Once a graph of transitions has been chosen, there still is
room for improvement. Indeed, one has some freedom in assigning
transition probabilities to each edge under the requirement, 
however, that the stationary distribution be of the right form.
It turns out that choosing appropriately those 
probabilities can lead to a sizable decrease in the mixing time.

In this context, Boyd et al.~\cite{BDX04} have recently observed 
that minimizing the mixing time of an MC on a graph of transitions
with a given stationary distribution
can be formulated as a semidefinite program (SDP), a well-known 
generalization of linear programming to matrices. See 
e.g.~\cite{BV03}. This enables the numerical computation
of the {\em fastest} mixing chain on a graph. Boyd et al.~\cite{BDX04}
have solved numerically a number of simple examples. 

A further benefit of this 
approach is that it provides a tight lower bound on
the optimal mixing time through the dual of the SDP.
In a follow-up paper, 
Boyd et al.~\cite{BDSX04} have used this bound
to exhibit an analytic expression for the fastest chain---and
prove its optimality---when 
the graph is made of a simple path under uniform distribution. 

However, a weakness of the SDP formulation is 
that only small graphs can be
studied thoroughly because numerical solvers run in time polynomial
in the size of the graph; in practice, chains have prohibitively
large state spaces. As for the dual, it is potentially useful from a
theoretical point of view even for complex chains, 
but Boyd et al.~\cite{BDX04} give no intuitive interpretation
of it, making it difficult to apply.

Our goal in this paper is to provide evidence 
that those shortcomings can be
overcome by a simpler approach. Our claim arises
from the following observation: one can
obtain lower and upper bounds on the mixing time of 
{\em completely specified} chains by way of well-known 
techniques such as path coupling, conductance, canonical 
paths etc.~\cite{J03}; formally, those bounds are
{\em parameterized by transition probabilities}. This prompts
the questions: can one optimize those bounds as functions
of the transition probabilities, and how close to optimum can one
get by doing so?

\subsection{Our results}\label{section:results}

We show through general results and examples 
that for well-structured problems, the above scheme can be
implemented, and that it is capable of providing nontrivial, sharp bounds.   

On the lower bound side, we use a standard extremal 
characterization to derive a general lower bound
which has a simple geometrical interpretation.
It consists in embedding the nodes of the graph into
an Euclidean space so as to stretch the nodes as much as
possible under constraints on the distance separating nodes
connected by an edge.
We show through convex optimization arguments that
it is actually tight. The simple interpretation makes it much easier
to apply than the dual SDP mentioned above.
Our result is similar to a bound obtained
recently
by Sun et al.~\cite{SBXD04} in a different
context. We also specialize the usual conductance bound
to the context of fastest mixing. We apply those
general results to several examples
obtaining close-to-optimal lower bounds.

On the upper bound side, it seems much harder to derive 
useful, general results.
A trivial bound can be obtained by considering any chain
on the graph, e.g. a canonical Metropolis-Hastings chain, 
and computing an upper bound on its mixing time. But as was shown by
Boyd et al.~\cite{BDX04}, there can be a large (unbounded) gap 
between standard and optimal chains. Instead, 
we show through examples that one can obtain
almost tight bounds by studying closely standard canonical paths
arguments and minimizing the bound over transition probabilities.
Put differently, our technique consists in identifying 
bottleneck edges and increasing the flow on them.
The fact that this scheme can work on nontrivial Markov chains
is not obvious {\sl a priori}, and this constitutes
our main result in the upper bound case. Moreover, this technique
is constructive and it allows to design a chain which might
be close to the fastest one. Our scheme is likely to work
only on well-structured problems but, even in that
case, there is no other non-numerical approach known---and
the numerical approach breaks down on large-scale problems.

Our main example is the Glauber dynamics of the Ising model, a problem
which is beyond the reach of the numerical SDP approach. 
In the case of the tree, by a judicious choice
of rates at which nodes are updated, we improve the mixing time
by an optimal factor.

\subsection{Organization of the paper}

We begin in Section~\ref{section:prelim} with a description of the setting
and approach of \cite{BDX04}. We introduce our main
techniques in Sections~\ref{section:lower} and~\ref{section:upper}. 
Section~\ref{section:rates} is devoted to optimal rates
of the Glauber dynamics of the Ising model.

\section{Preliminaries}\label{section:prelim}

\subsection{Setting}\label{section:rapid}

We are given an undirected graph $\gcal=(\vcal,\ecal)$
and a probability distribution
$\pi$ defined on the nodes of $\gcal$. We seek to sample from $\pi$ and
do so by running a reversible Markov chain $(X_t)_{t\geq 0}$ 
on the state space $\vcal$ with stationary distribution $\pi$, i.e. 
if $P = (P(i,j))_{i,j\in \vcal}$
denotes the transition matrix of $(X_t)_{t\geq 0}$, we must have
$\pi(i)P(i,j) = \pi(j)P(j,i),\ \forall i,j \in \vcal$.
We also require that the only transitions allowed are those given
by edges of $\gcal$, i.e. $P(i,j)=0,\ \forall (i,j)\notin \ecal$. 
For convenience, we assume that all self-loops are present.

The time to reach stationarity is 
governed by the second largest eigenvalue of $P$. 
More precisely, let $n = |\vcal|$ and
$1=\lambda_1(P)
> \lambda_2(P) \geq \cdots \geq \lambda_{n}(P) \geq -1$
be
the eigenvalues of $P$. We measure the
speed at which stationarity is reached by the {\sl relaxation time}
$\tau_2(P) = \frac{1}{1-\lambda_2(P)}$. See~\cite{AF04} for a thorough
discussion of other related quantities.
The smaller
$\lambda_{2}(P)$---and therefore $\tau_2(P)$---is, 
the faster $(X_t)$ approaches $\pi$. Given this
observation, it is natural to define the
{\sl fastest mixing chain} on $(\gcal, \pi)$ as the solution
of the optimization problem
\begin{equation}\label{eq:star}
\begin{array}{rlll}
& \multicolumn{3}{l}{\displaystyle\min_{P\geq 0}\ \lambda_{2}(P)}\\
&\mathrm{s.t.} & P(i,j) = 0, & \forall (i,j)\notin \ecal\\
&& \sum_{j\in \vcal} P(i,j) = 1, & \forall i \in \vcal\\
&& \pi(i)P(i,j) = \pi(j)P(j,i), & \forall (i,j)\in \ecal.
\end{array}
\end{equation}
In the remainder of this paper,
we save the notation $P^\star$ for a solution of (\ref{eq:star})---which 
might not be unique---and 
let $\lambda_2^\star = \lambda_2(P^{\star})$, and $\tau_2^\star 
= \tau_2(P^\star)$.
Note that our definition of fastest mixing differs
slightly from that in~\cite{BDX04}. Here, we take the usual approach
of ignoring the {\sl smallest} eigenvalue by considering the possibility
of adding a constant probability to each self-loop afterwards in order to 
bound the smallest eigenvalue away from $-1$.

\subsection{Fastest mixing via SDP}\label{section:fastest}

The main observation in~\cite{BDX04} is that (\ref{eq:star}) 
is actually
a semidefinite program (SDP). See e.g.~\cite{BV03} for background on
convex and semidefinite programming.
This observation makes possible the 
numerical computation of optimal transition
matrices. Unfortunately, since the running time of SDP algorithms
is at best polynomial in the size of the state space, this allows
only to study small graphs---for which sampling is actually
quite trivial. One idea put forward by Boyd et al.~\cite{BDX04} is to
solve the SDP on small instances of large combinatorial
problems and try and guess the structure of the optimal matrix
from the results. This is the approach used in~\cite{BDSX04}
to identify the optimal chain on the path. The prospect
of reproducing this type of exact result in other cases
seems limited.

From a theoretical point of view, an interesting consequence
of the SDP formulation is the existence of a dual which
can be used to give lower bounds on the optimal mixing time. 
Let $\|Y\|_*$ be the sum of the singular values of $Y$. Then,
in the case of the uniform stationary distribution,
the dual (of the more general version taking into account
the smallest eigenvalue) has the form~\cite{BDX04}
\begin{equation}\label{eq:minus}
\begin{array}{rlll}
& \multicolumn{3}{l}{\displaystyle\max_{z, Y}\ \sum_{i=1}^{n} z(i)}\\
&\mathrm{s.t.} & z(i) + z(j) \leq 2 Y(i,j), & \forall (i,j)\in \ecal\\
&& \sum_{j\in \vcal} Y(i,j) = 0, & \forall i \in \vcal\\
&& Y = Y^T,\quad \|Y\|_* \leq 1.
\end{array}
\end{equation}
Any feasible solution of (\ref{eq:minus}) provides a lower bound on the best
mixing time achievable on $(\gcal,\pi)$. Moreover, strong duality holds.
In~\cite{BDSX04}, 
this is used to prove optimality of a conjectured fastest chain
when the graph is a path. Note that giving an intuitive interpretation
of this optimization problem is not straightforward. 
This is a potential obstacle to the devising of good feasible solutions.

\section{Lower bounds}\label{section:lower}

In this section, we discuss general lower bounds on fastest mixing
that can be derived from common techniques for completely specified chains.
We apply our bounds to several examples.

\subsection{Variational characterization}\label{section:variational}

The standard lower bound for completely specified chains 
is based on a variational characterization
of the second eigenvalue of the transition matrix. See e.g.~\cite{AF04}.
To reveal the geometric flavor of our result, we will consider a
more general bound. Let $\psi_1, \ldots, \psi_n:\vcal \to \real$ be functions 
with 0 expectation 
under $\pi$, i.e.
$\sum_{i\in \vcal}\pi(i)\psi_l(i) = 0$ for all $l$ (where, as before, $n$ is
the number of nodes). For all $i\in\vcal$,
think of $\Psi(i) = (\psi_1(i),\ldots,\psi_n(i))$ as a vector associated to node $i$. 
Therefore,
$\Psi(1), \ldots, \Psi(n)$ is an embedding of the graph into $\real^n$.
For each $l$ separately, we have the inequality
\bea
(1 - \lambda_2(P))\sum_{k\in \vcal} \pi(k)\psi_l^2(k) 
\leq \sum_{(i,j)\in \ecal} (\psi_l(i) - \psi_l(j))^2 Q(i,j),  
\eea
where $Q(i,j) = \pi(i)P(i,j)$. Summing over $l$ we get the bound
\bea
1 - \lambda_2(P) 
\leq \frac{\sum_{(i,j)\in \ecal} \|\Psi(i) - \Psi(j)\|^2 Q(i,j)}
{\sum_{k\in \vcal} \pi(k)\|\Psi(k)\|^2},  
\eea
where $\|\,\cdot\,\|$ denotes the Euclidean norm in $\real^n$.
To turn the r.h.s. into
a bound on $1-\lambda_2^\star$,
we maximize over $Q$. But note that, for $\psi_1,\ldots,\psi_n$ 
fixed, the r.h.s. is linear
in $Q$ so this can be expressed as the {\sl linear} program
\begin{equation}\label{eq:primal}
\begin{array}{rclll}
1 - \lambda_2^\star &\leq& \multicolumn{3}{l}{\displaystyle\max_{Q\geq 0}\ 
\sum_{(i,j)\in \ecal} \frac{\|\Psi(i) - \Psi(j)\|^2}
{\sum_{k\in \vcal} \pi(k)\|\Psi(k)\|^2} Q(i,j)}\\
&&\mathrm{s.t.} & Q(i,j) = 0, & \forall (i,j)\notin \ecal\\
&&& \sum_{j\in \vcal} Q(i,j) = \pi(i), & \forall i \in \vcal\\
&&& Q(i,j) = Q(j,i), & \forall (i,j)\in \ecal.
\end{array}
\end{equation}
The dual of this linear program is\footnote{To obtain this particular form, one needs to consider
only those $Q(i,j)$'s such that $(i,j)\in \ecal$ and then only one of
$Q(i,j)$ and $Q(j,i)$.}
\begin{equation}\label{eq:linear}
\begin{array}{rclll}
1 - \lambda_2^\star &\leq& \multicolumn{3}{l}{\displaystyle\min_{\mathbf{z}}\ 
\sum_{i=1}^{n} \pi(i)z(i)}\\
&&\mathrm{s.t.} & z(i) + z(j) \geq \frac{\|\Psi(i) - \Psi(j)\|^2}
{\sum_{k\in \vcal} \pi(k)\|\Psi(k)\|^2}, & \forall (i,j)\in \ecal.
\end{array}
\end{equation}
Note the similarity with (\ref{eq:minus}). Note also that 
we can now minimize over $\psi_1,\ldots,\psi_n$ as well to get the best
bound possible. Make the change of variables 
$w(i) = z(i)\sum_{k\in \vcal} \pi(k)\|\Psi(k)\|^2$ for all $i\in\ecal$,
assume w.l.o.g. that $\sum_{i\in\vcal}\pi(i)w(i) = 1$ (one can always 
renormalize
the $\Psi$'s by $\sum_{i\in\vcal}\pi(i)w(i)$) and take the multiplicative
inverse of
the objective function. This finally leads to:
\begin{proposition}\label{prop:lowerbound}
The optimal relaxation time on $(\gcal,\pi)$ is bounded from below by
\begin{equation}\label{eq:lowerbound}
\begin{array}{rclll}
\tau_2^\star &\geq& \multicolumn{3}{l}{
\displaystyle\max_{\mathbf{w}, \Psi(1),\ldots,\Psi(n)}\ 
\sum_{k\in \vcal} \pi(k)\|\Psi(k)\|^2}\\
&&\mathrm{s.t.} &  \|\Psi(i) - \Psi(j)\|^2 
\leq w(i) + w(j), & \forall (i,j)\in \ecal\\
&&&\sum_{k=1}^n \pi(k)\Psi(k) = \mathbf{0}\\
&&&\sum_{i=1}^{n} \pi(i)w(i) = 1.
\end{array}
\end{equation}
Moreover, this bound is tight, i.e. we have equality above.
\end{proposition}
Informally, we seek to embed the graph into $\real^n$ so as to spread
the nodes as much as possible under local constraints over the
distances separating nodes connected by edges. The $w$'s give
some slack in choosing which edges are bound by stronger or
weaker constraints. See the examples. This bound is similar to that
obtained recently by~\cite{SBXD04} in a continuous-time context.
There, however, the r.h.s. in the inter-node distance
constraint is a {\sl fixed} weight $d_{ij}$ (instead of $w(i) + w(j)$), 
giving rise to a quite different problem.

\medskip

\noindent {\bf Proof (of tightness):} This follows from
convex optimization duality. To see this, we go back to formulation
(\ref{eq:linear}). Note that w.l.o.g., we can assume that $\sum_{k\in \vcal} \pi(k)\|\Psi(k)\|^2 = 1$. Make the change of variables $w(i) = \|\Psi(i)\|^2
- z(i)$ for all $i\in\vcal$, change the objective to 
$1 - \sum_{i=1}^{n} \pi(i)z(i) = \sum_{i=1}^{n} \pi(i)w(i)$, and 
set  $Y(i,j) = \Psi(i)^T\Psi(j)$ for all $i,j\in\vcal$. Then,
using
the Gram matrix representation for symmetric positive semidefinite matrices
(an $n\times n$ matrix $M$ is symmetric positive semidefinite 
if and only if there is a set
$\bfx_1, \ldots, \bfx_n$ of vectors in $\real^n$ such that 
$M_{ij} = \bfx_j^T\bfx_i$; see e.g.~\cite{HJ85}), we get
the equivalent bound
\begin{equation}\label{eq:linearagain}
\begin{array}{rclll}
\lambda_2^\star &\geq& \multicolumn{3}{l}{\displaystyle\min_{\mathbf{w}, Y=Y^T \succeq 0}\ 
\sum_{i=1}^{n} \pi(i)w(i)}\\
&&\mathrm{s.t.} & w(i) + w(j) \leq 2 Y(i,j),\quad \forall (i,j)\in \ecal &\\
&&& \sum_{i,j\in\vcal} \pi(i)\pi(j)Y(i,j) = 0,\ 
\sum_{k=1}^n \pi(k)Y(k,k) = 1,
\end{array}
\end{equation}
where $A \succeq 0$ indicates that $A$ is positive semidefinite. 
One can check that the dual of this convex optimization problem is
equivalent to minimizing the second largest eigenvalue over
reversible transition matrices on $(\gcal,\pi)$.
$\blacksquare$

Contrary to the standard setting, the multidimensionality
of the embedding seems necessary in the fastest mixing context. 
In particular, 
plugging the eigenvector corresponding to the second largest
eigenvalue of the optimal matrix as $\psi_1$ (with all
other coordinates 0) into (\ref{eq:lowerbound}) does not 
necessarily give a tight bound
because there is no guarantee that the optimal $w$'s will allow
enough room for a 1-dimensional embedding to spread sufficiently.

\begin{remark}
The above bound is actually very similar to that in the case
of completely specified chains which can be reformulated as
\begin{equation}\label{eq:specified}
\begin{array}{rclll}
\tau_2(P) &\geq& \multicolumn{3}{l}{
\displaystyle\max_{\Psi(1),\ldots,\Psi(n)}\ 
\sum_{k\in \vcal} \pi(k)\|\Psi(k)\|^2}\\
&&\mathrm{s.t.} &  \sum_{(i,j)\in\ecal}\|\Psi(i) - \Psi(j)\|^2Q(i,j) = 1\\
&&&\sum_{k=1}^n \pi(k)\Psi(k) = \mathbf{0}.
\end{array}
\end{equation}
Here the ``slack'' takes the form of a fixed weighted average over inter-node
distances. The multidimensionality turns out not to be necessary in this case.
\end{remark}
\begin{remark}
The same scheme can be applied to the log-Sobolev constant. In that
case, one maximizes the entropy instead of the variance. See also \cite{BDX04}.
\end{remark}
\begin{remark}
The smallest eigenvalue has its own geometry. There, the bound is
the same with the term
$\|\Psi(i) - \Psi(j)\|^2$ in the inter-node distance constraint
replaced by $\|\Psi(i) + \Psi(j)\|^2$. The formulation (\ref{eq:minus})
is equivalent to a combination of the two geometries (smallest and second largest 
eigenvalues).
\end{remark}

\subsection{Conductance}\label{section:conductance}

As an illustration of Proposition~\ref{prop:lowerbound} we give a simple 
adaptation of the conductance bound to the context of fastest mixing.
\begin{proposition}\label{prop:conductance}
Let $\Upsilon$ be the {\sl weighted vertex expansion} of $(\gcal,\pi)$ 
\bea
\Upsilon 
&\equiv& \min \left\{\frac{\pi(\delta S)}
{\pi(S)\land \pi(S^c)}
\ :\ S\subseteq \vcal \right\},
\eea
where $a\land b = \min\{a,b\}$ and 
$\delta S$ is the set of nodes $i\in S^c$ such that there is
a $j\in S$ with $(i,j)\in \ecal$. We have the following bound
\bea
\tau_{2}^{\star} 
&\geq& \frac{1}{2\Upsilon}.
\eea
\end{proposition}

\medskip

\noindent This
bound is actually folklore. It is easily derived from the usual conductance bound
and is often used to obtain lower bounds on completely specificied chains. Here we give
a direct proof.

\medskip

\noindent {\bf Proof:} A simple embedding of $\gcal$ in $\real^n$
is to map each node to one of only 2 points $\bfx_0,\bfx_1$. 
Say the subset $S \subseteq \vcal$
is mapped to $\bfx_0$. Then we must have $\pi(S)\bfx_0 + \pi(S^c)\bfx_1 = \mathbf{0}$.
Also, since the distance between nodes inside $S$ (resp. $S^c$) is 0, we can
set w.l.o.g. the $w$'s of nodes {\sl not on the boundary} of $S$ (resp.
$S^c$) to 0. We assign to the points on the boundary of $S$ (resp. $S^c$)
the value $w_0$ (resp. $w_1$). Since we care only about the sum $w_0 + w_1$
and the only constraint on $w_0, w_1$ is $\pi(\delta S^c)w_0 + \pi(\delta S)w_1 = 1$,
it is in our advantage to fix one of $w_0, w_1$ to 0 as well. 
Say $\pi(\delta S^c)\leq \pi(\delta S)$
w.l.o.g. Then $\|\,\bfx_0 - \bfx_1\|^2 = w_0 = (\pi(\delta S^c))^{-1}$ and $w_1 = 0$. An easy
calculation gives $\|\,\bfx_0\|^2 = (\pi(\delta S^c))^{-1}(1 - \pi(S))^2$
and $\|\,\bfx_1\|^2 = (\pi(\delta S^c))^{-1}(\pi(S))^2$. Therefore,
$\pi(S)\|\,\bfx_0\|^2 + \pi(S^c)\|\,\bfx_1\|^2 = \frac{\pi(S)\pi(S^c)}{\pi(\delta S^c)}$
and the result follows.
$\blacksquare$

\subsection{Examples}

\subsubsection{$K_n - K_n$}

This is the graph made of two $n$-node complete graphs
joined by an edge. 
We denote the nodes on one side of the
linking edge by
$1,\ldots,n$ and those on the other side by $1',\ldots,n'$. The linking
edge is $(1,1')$.
The stationary distribution
is uniform. The vertex expansion bound gives
$\Upsilon = \frac{1/(2n)}{n/(2n)} = \frac{1}{n}$ and
$\tau_2^\star \geq n/2$.
To get something sharper, 
we appeal to our more
general bound. The bottleneck in this graph is intrinsically one-dimensional, so
we take all coordinates except the first one to be 0, i.e. we consider only $\psi_1$.
By symmetry, it is natural to map the nodes to $\psi_1(1) = -\psi_1(1')
= x_0$ and $\psi_1(i) = - \psi_1(i') = x_1$, for $i\neq 1$, 
with $0 \leq x_0 \leq x_1$. The main insight here is that we should make the
distance between $1$ and $1'$ as large as possible because that pushes away from 0
all the other points at the same time (because of the local constraints). 
So we take $w(i) = w(i') = 0$, for all
$i\neq 1$, and $w(1) = w(1') = n$, which gives $x_0 = \frac{\sqrt{2}}{2} \sqrt{n}$
and $x_1 = \frac{\sqrt{2} + 2}{2} \sqrt{n}$. Summing the squares leads to a lower
bound asymptotic to $(\frac{3}{2} + \sqrt{2})n \geq 2.914 n$. In 
Section~\ref{section:upper}, we give an almost matching upper bound. See 
also~\cite{BDPX04} for a similar upper bound.

\subsubsection{$n$-cycle and $d$-dimensional torus}

In constrast to our preceding example, the $n$-cycle gives rise
naturally to a multidimensional embedding. We let the stationary
distribution be uniform. By symmetry we choose
all $w$'s equal. So all pairs of consecutive nodes have to be embedded
to points at distance (at most) $\sqrt{2}$. Our goal of maximizing 
the sum of the squared 
norms---and the natural symmetry--- 
leads to spreading the points evenly on a circle centered around the origin
(in any 2-dimensional subspace of $\real^n$). That is, we take all
coordinates except the first two to be 0 and, numbering the nodes
from 1 to $n$ in order of traversal,
we let $(\psi_1(i), \psi_2(i)) = (R\cos(2\pi i/n),
R\sin(2\pi i/n))$, $i=1,\ldots,n$, 
for a value of $R$ which remains to be determined. 
The distance between consecutive points has to be $\sqrt{2}$ so a little
geometry suggests $R = \frac{\sqrt{2}}{2 \sin(\pi/n)} \geq \frac{\sqrt{2}n}{2\pi}$.
Thus the lower bound is $\tau_2^\star \geq \frac{n^2}{2\pi^2}$, matching
the relaxation time of the symmetric walk. See e.g.~\cite{AF04}.

One can generalize this result to the $m^d$-point grid on a $d$-dimensional
torus by considering a $2d$-dimensional embedding. For $1\leq i_1,\ldots,i_d \leq m$,
node $(i_1,\ldots,i_d)$
is mapped to
\bea
(R\cos(2\pi i_1/m),
R\sin(2\pi i_1/m), \ldots, R\cos(2\pi i_d/m),
R\sin(2\pi i_d/m)),
\eea 
with $R$ as above. Thus, $\tau_2^\star \geq \frac{dm^2}{2\pi^2}$, again matching
the relaxation time of the symmetric walk. See~\cite{AF04}.

\subsubsection{Geometric random graphs}

In their analysis of random walks on geometric random graphs,
Boyd et al.~\cite{BGPS04} consider, in a key step, a variant of the $d$-dimensional
grid of the previous example. Let $k$ be a fixed integer smaller than $m$.
Again, our graph is made of the $m^d$
points of the $d$-dimensional torus $\intgr_m^d$ (integers modulo $m$) with 
uniform stationary distribution.
Two nodes $(i_1,\ldots,i_d)$ and $(j_1,\ldots,j_d)$ are connected by an edge 
if $i_l - j_l$  modulo $m$ 
is less or equal to $k$ for all $1\leq l\leq d$ (the points are at most $k$
cells apart in every dimension). Because of the ``diagonal'' edges, it seems natural to collapse
all nodes on a single $m$-cycle. More precisely, we map $(i_1,\ldots, i_d)$
to $(R\cos(2\pi i_1/m),
R\sin(2\pi i_1/m))$. We take uniform $w$'s.
Because some edges connect nodes $k$ steps apart,
the radius (which is constrained by the fact that points connected by an edge are
at most $\sqrt{2}$ apart) is now $R = \frac{\sqrt{2}}{2 \sin(k\pi/m)} \geq \frac{\sqrt{2}m}{2k\pi}$ (assume that $k$ divides $n$ for convenience).
Thus $\tau_2^\star \geq \frac{m^2}{2k^2\pi^2} = \Theta((n/D_d)^{2/d})$,
where $D_d$ is the degree of each node and $n$ is the number of nodes.
This bound matches the lower bound in~\cite{BGPS04}. There, exact expressions
for the eigenvalues of tensor products of circulant matrices 
and the analysis of a linear program lead
to a lower bound on fastest mixing on this graph. Our geometric method is
much simpler.

\begin{remark}
In the previous two examples, plugging the same embeddings into the
completely specifed setting (\ref{eq:specified}) gives tight lower bounds
on the symmetric walks. More generally, the lower bound in 
Proposition~\ref{prop:lowerbound} applies to any completely
specified chain---as do all lower bounds on fastest mixing---and 
it could prove useful as an alternative
to the standard variational characterization when the precise
details of the transition matrix appear too
cumbersome.
\end{remark}

\section{Upper bounds}\label{section:upper}

It seems difficult to give general upper bounds on fastest mixing.
An obvious technique is to pick
an arbitrary chain and compute an upper bound on its relaxation time.
For example, one might use the canonical (max-degree like) chain
defined by the transition probabilities 
$P_d(i,j) = \pi(j)/\pi_*$  if $(i,j) \in \ecal$ (and
0 otherwise) with 
$\pi_* = \max\{\sum_{j:(i,j)\in\ecal} \pi(j)\,:\,i\in\vcal\}$.
Let $\pi_0 = \min_{i\in\vcal}\pi(i)$ and 
recall the definition of vertex expansion $\Upsilon$ from
Proposition~\ref{prop:conductance}. 
Noting that for any subset $S\subseteq\vcal$,
\bea
\sum_{i\in S}\sum_{j\in S^c} \pi(i)P_d(i,j)
\geq \left(\frac{\pi_0}{\pi_*}\right)\sum_{i\in \delta S^c}\sum_{j\in \delta S} 
\pi(i)\ind_{(i,j)\in\ecal} \geq \left(\frac{\pi_0}{\pi_*}\right) \pi(\delta S^c), 
\eea
and applying the standard Cheeger inequality to $P_d$ leads to,
\bea
\tau_2^\star \leq \left(\frac{\pi_*}{\pi_0}\right)^2\frac{2}{\Upsilon^2}.
\eea
A different chain would have provided a different---and possibly better---bound.
Anyhow, this Cheeger-type bound is very unlikely to lead to
useful results, and moreover it tells us nothing about the optimal chain. 

Instead, the goal of this section is to illustrate the computation
of a nontrivial upper bound through a canonical paths argument. 
The underlying idea is similar to that used in the lower bound
above. That is, we think of a standard upper bound for completely
specified chains as parameterized by transition probabilities
and attempt to minimize the bound over those probabilities.
It turns out that because
of its straightforward dependence on the transition matrix, the canonical
paths bound appears to be the most manageable. 
In this section and the next one,
we show by way of examples that it can actually lead to sharp results. 

\subsection{Canonical paths: $K_n - K_n$ example continued}\label{section:paths}

We consider again the $K_n-K_n$ graph with uniform distribution.
This chain is analyzed in details in~\cite{BDPX04},
where using sophisticated group-theoretic-based symmetry analysis, 
all eigenvalues are computed. Here, we give a very different, much more
elementary, treatment. Also, being simpler, 
our approach has the potential of being applicable
more generally.
We proceed as follows: we write down the canonical paths upper bound
as a function of $P$; we then choose $P$ among $\pi$-reversible chains
so as to minimize the bound. Given a set $\Gamma$ of paths $\gamma_{xy}$
in $\gcal$
for all pairs of nodes $x,y$, the canonical paths upper bound is
\bea
\tau_2(P) \leq \bar\rho(P, \Gamma),
\eea
with
\bean\label{eq:congestion}
\bar\rho(P, \Gamma)
&=& \max_{e} \frac{\sum_{\gamma_{xy} \ni e}\pi(x)\pi(y)|\gamma_{xy}|}{Q(e)},
\eean
where $|\gamma_{xy}|$ is the number of edges in $\gamma_{xy}$.
Notice that the choice of paths depends---crucially---only on the graph and is
therefore valid for any transition matrix consistent with $(\gcal,\pi)$.
Let $W(e)$ be the numerator in (\ref{eq:congestion}).
On $K_n-K_n$, the natural choice of paths is to let $\gamma_{xy}$
be the shortest path (in terms of number of edges) between $x$ and $y$.
Then
\bea
&& W(i,j) 
= \frac{1}{(2n)^2},\quad \forall i,j \neq 1,\\
&&W(i,1) 
= \frac{1}{(2n)^2}(1 + 2 + 3(n-1)) = \frac{3n}{(2n)^2}, 
\quad \forall i,\\
&& W(1,1')
= \frac{1}{(2n)^2}(1 + 4(n-1) + 3(n-1)^2)
= \frac{3n^2 - 2n}{(2n)^2}.
\eea
Similar values hold for the other complete subgraph. The largest
contribution to the maximum above clearly comes from $W(1,1')$.
In order to decrease the ratio in $\bar\rho(P, \Gamma)$,
we need to
choose a large value for $Q(1,1')$.
But as we increase $Q(1,1')$, the $Q(1,i)$'s
and $Q(1',i')$'s have to be lowered accordingly. 
We do so until congestion is the same
on edges $(1,1')$, $(1,i)$'s and $(1',i')$'s. That is, we require
\bea
&& \frac{W(1,1')}{Q(1,1')} = \frac{W(1,i)}{Q(1,i)},\quad \forall i\neq 1 \qquad
Q(1,1') + \sum_{i=2}^{n} Q(1,i) = \frac{1}{2n}, 
\eea
and similarly for the other side. The solution is
\bea
&& Q(1,i) = Q(1',i')
= \frac{1}{2n(2n- 5/3)},\quad \forall i\neq 1,\\
&& Q(1,1') 
= \frac{n - 2/3}{2n(2n - 5/3)}.
\eea
We extend this to all edges by
\bea
&& Q(i,j) = Q(i',j')
= \frac{1}{n-1}\left[1 - \frac{1}{2n(2n- 5/3)}\right],
\quad \forall i,j\neq 1.
\eea
The upper bound becomes
\bea
&&\tau_{2}^{\star} \leq 3n(1 - 5/(6n)).
\eea
Recall that our lower bound was $\tau_2^\star \geq 2.9n$.
Note that the standard chain would have consisted in choosing 
a neighbour
uniformly at random at each step. The same calulation gives an upper
bound of $\Omega(n^2)$ in that case.

\begin{remark}
In summary, our upper bound technique consists in two steps:
identify transitions contributing to slow mixing by computing
the congestion ratio in (\ref{eq:congestion}); then increase
as much as possible the probability of transition on those
bottleneck edges. Instead, one might try to use the same idea
with conductance (or other upper bounds). 
But in that case, the fact that all cuts---instead of edges---have
to be accounted for simultaneously makes the task more difficult.
\end{remark}

\section{Optimal rates for Glauber dynamics}\label{section:rates}

In this section, we show that the framework discussed so far
can be applied to large, well-structured combinatorial problems
where the numerical SDP method has little chance of being helpful. 

\subsection{Glauber dynamics}\label{section:glauber}

Let $G = (V, E)$ be a finite graph\footnote{We now have two graphs. As before,
calligraphic letters are used to denote the transition graph (see below).}. 
A configuration on $G$ is a map $\sigma : V \to C$,
where $C$ is a finite set. Typically, $\sigma$
is a spin or a color. We consider the following stationary distribution on 
$C^V$ 
\bea
\pi(\sigma)
&=& \frac{1}{Z}\prod_{(v,w) \in E} \alpha_{vw}(\sigma(v), \sigma(w)),  
\eea 
where $Z$ is a normalization constant and $(v,w)$ is an undirected
edge with endpoints $v,w$. Let $\scal \subseteq C^V$ be the subset of $C^V$
on which $\pi$ is nonzero. We wish to sample from $\pi$ by running
a reversible MC on $\scal$, but allow
only transitions that change the state of one node at a time, i.e. the
transition graph is
$\gcal = (\vcal, \ecal)$ with $\vcal = \scal$ and $(\sigma, \sigma') \in \ecal$ 
if and only if
$\sigma(v) = \sigma'(v)$ for all but at most one node $v \in V$. Let
$\sigma_v^a$ be the configuration
\begin{equation*}
\sigma_v^a(w) 
= \left\{\begin{array}{ll}
\sigma(w), & \mathrm{if}\ w\neq v,\\
a,& \mathrm{if}\ w = v.
\end{array}
\right.
\end{equation*}
One such ``local'' MC is the so-called Glauber dynamics which, at each step, picks
a node $v$ of $G$ uniformly at random and updates the value
$\sigma(v)$ according to the transition probability distribution
\bea
K(\sigma,\sigma_v^a)
&=& \frac{\displaystyle\prod_{w:(v,w)\in E} \alpha_{vw}(a,\sigma(w))}
{\displaystyle\sum_{a'\in C} \prod_{w:(v,w)\in E} \alpha_{vw}(a',\sigma(w))}.
\eea 
One can check that $K$ is $\pi$-reversible. We actually consider a
generalization of the Glauber dynamics by allowing the update rates to vary. 
More precisely, at each step, we pick a node
$v$ of $G$ with probability $\rho(v)$ for some distribution
$\rho : V \to [0,1]$, and we
update $\sigma(v)$ according to $K$ as above. The standard chain corresponds
to uniform $\rho$. 

Predictably the question we ask is: can we compute the
rates $\rho$ minimizing the mixing time? Or at least can we get reasonable lower and
upper bounds on fastest mixing in this restricted setting? 
We do so by following the methodology put forward
in the previous sections.

We first give an elementary bound on the best achievable improvement. This
observation is essentially due to~\cite{BDX04}.
\begin{proposition}\label{prop:notbad}
Let $P^\star$ be the fastest chain on $(\gcal, \pi)$ (not necessarily
of the Glauber dynamics type). Also, let $\rho^\star$ (resp. $U$) 
be the optimal (resp. uniform) rates
for the Glauber dynamics. Denote by $P_\rho$ the Glauber dynamics with 
rates $\rho$ and let
$\overline{K} = \max_{\sigma, v, a} K^{-1}(\sigma,\sigma_v^a)$.
Then,
\bea
\tau_2(P_U)
&\leq& \overline{K}\, |V| \tau_2(P^\star),\qquad
\mathit{and}\qquad 
\tau_2(P_U) \leq |V|\tau_2(P_{\rho^\star}).
\eea
\end{proposition}
\noindent{\bf Proof:} By the variational
characterization of $\lambda_2(P^\star)$ and the fact that
$P^\star(\sigma, \sigma') \leq 1$, 
\bea
1 - \lambda_2(P^\star)
&=& \inf_g \frac{\displaystyle\sum_{\sigma} \sum_{v\in V}\sum_{a\in C}
(g(\sigma) - g(\sigma_v^a))^2 \pi(\sigma) P^\star(\sigma, \sigma_v^a)}
{\displaystyle\sum_{\sigma} \sum_{v\in V}\sum_{a\in C}
(g(\sigma) - g(\sigma_v^a))^2 \pi(\sigma) \pi(\sigma_v^a)}\\
&\leq& \overline{K}\, |V| \  
\inf_g \frac{\displaystyle\sum_{\sigma} \sum_{v\in V}\sum_{a\in C}
(g(\sigma) - g(\sigma_v^a))^2 \pi(\sigma) U(v) K(\sigma, \sigma_v^a)}
{\displaystyle\sum_{\sigma} \sum_{v\in V}\sum_{a\in C}
(g(\sigma) - g(\sigma_v^a))^2 \pi(\sigma) \pi(\sigma_v^a)}\\
&=& \overline{K}\, |V| \  (1 - \lambda_2(P_U)).
\eea
A similar argument
gives the second inequality.
$\blacksquare$

\noindent Thus, assume $\overline{K}$ is $O(1)$, then
the best improvement over $P_U$
one can hope for is a factor of $O(|V|)$.

We now use a canonical paths argument similar to that in
Section~\ref{section:upper} to obtain a general upper bound
on fastest mixing for Glauber dynamics. 
\begin{proposition}\label{prop:rates}
Let $\Gamma$ be a set of paths $\gamma_{\sigma,\sigma'}$ 
in $\gcal$ for each pair $\sigma,\sigma'$
in $\scal$.  
Assume we have a bound $B_v$ 
(depending only on $v$) on the ratio appearing
in the canonical paths bound (\ref{eq:congestion}) 
for edges of the form $(\sigma,\sigma_v^a)$ in the uniform rates
case. Then,
\bea
\tau_2(P_U) \leq \max_v B_v,\qquad
\mathit{and}\qquad
\tau_2(P_{\tilde\rho}) \leq |V|^{-1}\sum_v B_v,
\eea
with the choice of rates
$\tilde\rho(v) 
= B_v / \sum_u B_u$.
\end{proposition}
\noindent{\bf Proof:} The first inequality is the canonical paths
bound. For the second one, note that the ratio
in (\ref{eq:congestion}) is multiplied by $(|V| B_v/\sum_u B_u)^{-1}$
when replacing uniform rates with $\tilde\rho(v)$. We then apply
the canonical paths bound to $P_{\tilde\rho}$ using the bound
$B_v$ and the previous observation. Note that $\tilde\rho$
is the choice of rates that makes all bounds on the ratio in
(\ref{eq:congestion}) equal.
$\blacksquare$

The point
of Proposition~\ref{prop:rates} is that optimal improvement can be attained
if most $B_v$'s are small compared to $\max_v B_v$. We give such an example
in the next subsection.

\subsection{Special case: the Ising model}\label{section:ising}

We apply the previous result to the case of the Ising model on a finite graph. 
Here $C = \{-1,+1\}$, $\scal = C^V$, and
$\alpha_{vw}(\sigma(v), \sigma(w)) = \exp\left(\beta \sigma(v)\sigma(w)\right)$,
where $\beta > 0$ is some constant.

As shown in~\cite{KMP01}, the mixing time of the Glauber dynamics on
a graph depends on its cut-width.
\begin{definition}
The cut-width $\xi(G)$ of a graph $G$ is the smallest integer
such that there exists a labeling $v_1,\ldots,v_{|V|}$ of the vertices
such that for all $1\leq k\leq |V|$ the number of edges from
$\{v_1,\ldots,v_{k}\}$ to $\{v_{k+1},\ldots,v_{|V|}\}$ is at most
$\xi(G)$.
\end{definition}
To use Proposition~\ref{prop:rates}, we have to define the width of each
node. Let $I : V \to \{1,\ldots,|V|\}$ be some ordering of the nodes 
(not necessarily optimal), then we let $\xi^I(v)$ be the number of edges
from $\{w\ :\ I(w) \leq I(v) \}$ to $\{w\ :\ I(w) > I(v)\}$.
Let $\Delta$ be the maximum degree of $G$. Then 
it follows from~\cite{KMP01} that
a bound as required in Proposition~\ref{prop:rates} is 
\bea
B_v \equiv |V|^2\,\exp\left((4 \xi^I(v) + 2\Delta)\beta\right),
\eea
with in particular 
$\max_v B_v = |V|^2\, \exp\left((4 \xi(G) + 2\Delta)\beta\right)$ if
$I$ is an optimal ordering.

One can try and compute $\sum_v B_v/|V|$ in special cases. A 
rather uninteresting
graph is the $s\times s$ grid. There, a natural ordering is to 
start from a corner, move horizontally as far as one can, then go to the next line
and start over. In this ordering, the width of most nodes, including the 
maximum-width node, is approximately $s$ and 
therefore using non-uniform rates has essentially no effect.   

Here is a more interesting example. 
Let $T_r^{(b)} = (V_r, E_r)$ be the complete rooted $b$-ary tree with $r$ levels
(the root is at level $0$ and the leafs, at level $r$). 
Let $n_r$ be the number of vertices in $T_r^{(b)}$.
\begin{proposition}\label{prop:tree}
For $\beta$ large enough, an appropriate choice of rates 
leads to the estimate
\bea
\tau_2(P_{\rho^\star}) =
O\left(n_r\,e^{4(b-1)\beta r}\right),
\eea
as $r$ tends to $+\infty$.
In constrast, the best known upper bound on the uniform Glauber dynamics~\cite{KMP01}
is
\bea
\tau_2(P_{U}) =
O\left(n_r^{2}\,e^{4(b-1)\beta r}\right).
\eea
\end{proposition}
\noindent {\bf Proof:}
A good ordering of nodes of $T_r^{(b)}$, say $I$, is given by a depth-first search (DFS) 
traversal of the tree starting from the root.
This implies that $\xi(T_r^{(b)}) < (b-1)r + 1$~\cite{KMP01}.
Note that the width of a node $v$ is the number of unvisited neighbours 
of previously visited vertices when the DFS search reaches $v$. Therefore,
the width of the root
is $b$. Then, say vertex $v$ is on level $1 \leq l < r$ and is
the $q$-th child of its parent $w$ (in the DFS traversal order). Then 
$\xi^I(v) = \xi^I(w) + b - q$ because (1) $v$ has $b$ children, (2)
$q$ children of $w$ have now been visited, and (3) all descendants of the
first $q-1$ children of $w$ have been visited---so these add nothing
to the width. As for nodes on level $r$, we have similarly
$\xi^I(v) = \xi^I(w) - q$ if $v$ is the $q$-th child of $w$. Thus,
the contribution to $\sum_v B_v$ of the $l$-th level, $1 \leq l < r$,
is 
\bea
B^{(l)} 
&=& B^{(l-1)}\left(e^{4(b-1)\beta} + \cdots + e^{4(0)\beta}\right)\\
&=& B^{(l-1)}\left(\frac{e^{4 b \beta} - 1}{e^{4 \beta} - 1}\right)
\equiv B^{(l-1)} \zeta(b, \beta), 
\eea
with a similar expression for $l=r$. Summing over all levels, we get
\bea
\frac{\sum_v B_v}{|V|}
&=& n_r^2\, \frac{e^{(4b + 2\Delta)\beta}}{n_r}\,
\left(1 + \zeta(b,\beta) + \cdots + \zeta(b,\beta)^{r-1}
+ e^{-4b\beta}\zeta(b,\beta)^r\right)\\
&=& n_r\, e^{(4b + 2\Delta)\beta}\, 
\left(\frac{\zeta(b,\beta)^r - 1}{\zeta(b,\beta) - 1}
+ e^{-4b\beta}\zeta(b,\beta)^r\right).
\eea
In the low-temperature regime, i.e. for $\beta$ large
(we actually assume $e^{4\beta} \gg 1$), this is
\bea
\frac{\sum_v B_v}{|V|}
&=& O\Big(n_r\, \exp\Big\{4[(b-1)r + 1]\beta + 2\Delta\beta\Big\}\Big),
\eea 
whereas 
\bea
\max_v B_v &=& 
n_r^2\, \exp\left((4 \xi(T_r^{(b)}) + 2\Delta)\beta\right)\\ 
&=&
n_r^2\, \exp\Big\{4[(b-1)r + 1]\beta + 2\Delta\beta\Big\}.
\eea
Therefore, we get an optimal improvement of $O(n_r)$ over the usual
Glauber dynamics. 
$\blacksquare$

For a lower bound, we have the following result where
we assume $b=3$ for convenience.
\begin{proposition}\label{prop:lowertree}
Assume $b=3$ and let $\epsilon = (1 + e^{2\beta})^{-1}$. Then
\bea
\tau_2^\star =
\Omega\left(n_r^{-\ln(2\epsilon + 8\epsilon^2) - 1}\right),
\eea
as $r$ tends to $+\infty$.
In constrast, the best known lower bound in the uniform
case~\cite{KMP01} is
\bea
\tau_2(P_U) =
\Omega\left(n_r^{-\ln(2\epsilon + 8\epsilon^2)}\right).
\eea
\end{proposition}
\noindent {\bf Proof:}
Kenyon et al.~\cite{KMP01} use 
recursive majority to define a cut in the space of configurations
and apply the conductance bound.
The recursive majority $m(\sigma)$
of a configuration $\sigma$ is computed as follows: set
$M(v) = \sigma(v)$ for all $v$ on level $r$; starting from level
$r-1$ and up, compute
$M$ on each node by taking the majority
of the values of $M$ at the children of that node; 
output the value of $M$ at the root.
Let $S$ be the set
of configurations $\sigma$ with $m(\sigma) = +1$.  It follows
from~\cite{KMP01} that, under $\pi$, the probability that
a configuration is such that its recursive majority
is flipped by changing the value at a fixed leaf is
at most $(2\epsilon + 8\epsilon^2)^{r-1}$.
The union bound and the $\{-1,+1\}$ symmetry
imply that
$\pi(\delta S^c) \leq \frac{3^r}{2}(2\epsilon + 8\epsilon^2)^{r-1}$ and
$\pi(S) = \frac{1}{2}$.  
By Proposition~\ref{prop:conductance}, we deduce
$\lambda_{2}^\star \geq 1 - 2 (3)^r(2\epsilon + 8\epsilon^2)^{r-1}$.
On the other hand, the usual conductance bound applied to 
the uniform case gives that 
$\lambda_{2}(P_U) \geq 1 - 2 \Phi_S$, with
\bea
\Phi_S 
&=& \pi(S)^{-1}\displaystyle\sum_{\sigma\in S, \tau \in S^c} \pi(\sigma) P_U(\sigma,\tau)
\leq 2(3)^{-r}\displaystyle\sum_{\genfrac{}{}{0pt}{}{\sigma\in S, \tau \in S^c}
{(\sigma,\tau)\in \ecal}} \pi(\sigma)
\leq (2\epsilon + 8\epsilon^2)^{r-1},
\eea
where we have used that $P_U(\sigma,\tau) \leq 3^{-r}$ for
neighbours $\sigma,\tau$~\cite{KMP01}. Since $3^r = O(n_r)$ 
our lower bound on fastest mixing is $O(n_r)$ times smaller
than that on the standard Glauber dynamics. 
$\blacksquare$

Obtaining tighter bounds would require a sharper analysis in the standard
setting.

\begin{remark}
We are not claiming that this choice of rates leads to the
fastest sampling algorithm for this model. Indeed, in the case
of the Ising model on a tree, a very simple propagation algorithm
is much faster~\cite{EKPS00}.
Rather, our point is to establish that fastest mixing analysis
is feasible on nontrivial large-scale chains---a fact that was not
immediate from previous works. It remains to be seen whether 
fastest mixing ideas will find
useful applications in sampling. 
\end{remark}

\section*{Acknowledgements}

This work was motivated by a talk of Persi Diaconis. 
We thank David Aldous, Elchanan Mossel, and Santosh Vempala for comments.
Part of this work was done while visiting CSAIL at MIT. The author
acknowledges the partial support of NSERC.

\end{document}